\documentclass[12pt]{amsart}

\usepackage{amsfonts}
\usepackage{amsthm}     % need for theorem environments
\usepackage{amsmath}    % need for subequations
\usepackage[dvips]{graphicx}   % need for figures
\usepackage{verbatim}   % useful for program listings
\usepackage{color}      % use if color is used in text
\usepackage{subfigure}  % use for side-by-side figures
\usepackage{hyperref}   % use for hypertext links, including those to external documents and URLs

\newtheorem{thm}{Theorem}[section]

\newtheorem{lemma}[thm]{Lemma}

\theoremstyle{definition}
\newtheorem{defn}[thm]{Definition}

\theoremstyle{remark}

\numberwithin{equation}{section}

\begin{document}
\title[Complementary Regions of Knot and Link Diagrams]
{Complementary Regions of Knot and Link Diagrams}

\date{\today}
\author[Colin Adams]{Colin Adams}
\address{Colin Adams, 
Department of Mathematics and Statistics, Williams 
College, Williamstown, MA 01267}
\email{Colin.C.Adams@williams.edu}

\author[Reiko Shinjo]{Reiko Shinjo}
\address{Reiko Shinjo, 
Osaka City University Advanced Mathematical Institute, 
3-3-138 Sugimoto, Sumiyoshi-ku, Osaka 558-8585, Japan}
\email{reiko@suou.waseda.jp}

\author[Kokoro Tanaka]{Kokoro Tanaka}
\address{Kokoro Tanaka, 
Department of Mathematics, Tokyo Gakugei University, 
Nukuikita 4-1-1, Koganei, Tokyo 184-8501, Japan}
\email{kotanaka@u-gakugei.ac.jp}
\maketitle

\begin{abstract}
An increasing sequence of integers is said to be universal for knots  
and links if every knot and link has a projection to the sphere such  
that the number of edges of each complementary face of the projection  
comes from the given sequence. In this paper, it is proved that the  
following infinite sequences are all universal for knots and links:  
$(3,5,7, \dots)$, $(2,n, n+1,n+2,\dots)$ for all $n \geq 3$ and 
$(3,n,n+1,n+2, \dots)$ for all $n \geq 4$. 
Moreover, the following finite sequences are also universal 
for knots and links: $(3,4,5)$ and $(2,4,5)$.  
%Moreover, the finite sequence $(3,4,n)$, for all $n\geq 5$, 
%is universal for knots, and the following finite sequences 
%are also universal for knots and links: $(3,4,5)$ and $(2,4,5)$.
%
It is also shown that every knot has a projection with exactly 
two odd-sided faces, which can be taken to be triangles, 
and every link of $n$ components has a projection with at  
most $n$ odd-sided faces if $n$ is even and $n+1$ odd-sided faces 
if $n$ is odd.
\end{abstract}

\section{Introduction}\label{S:intro}

Knot diagrams have proved to be an extremely powerful tool for analyzing knots and links in 3-space. They can be used to tabulate knots and links as in \cite{HTW}. They can be used to define and compute invariants, such as the linking number or the various polynomials associated to knots. But relatively little seems to be known about the combinatorial properties of the knot diagrams associated to knots.

Given a diagram of a knot  or link $K$, one can ignore which strand is the overstrand at each crossing and think of it as a planar 4-valent graph embedded on the 2-sphere $S^2$. (For the purposes of this paper, links refer only to multi-component links and do not include knots.) This graph divides the sphere into $n$-gons, which we call faces and which meet at vertices and along edges. For our purposes, we  assume that the diagrams are connected and reduced, which for the diagram means that any obviously unnecessary crossings have been eliminated, and which for the graph means that the removal of any single vertex cannot separate the graph. In particular, the two endpoints of an edge cannot both occur at the same vertex.

 One would like to understand the possibilities for the collection of complementary $n$-gon faces associated to the diagrams of a knot or link. To that end, we make the following definitions.

\begin{defn} Given a knot or link $K$ and a strictly increasing sequence of integers $(a_1, a_2, a_3, ...)$ with $a_1 \geq 2$, we say  the sequence is {\it realized} by $K$ if there exists a diagram for the knot or link such that each face is an $a_n$-gon for some $a_n$ that appears in the sequence. Such a diagram is called an $(a_1, a_2, a_3, ...)$-diagram. (Not every $a_n$ must be realized by a face.) We say that a sequence is  {\it universal}  if every knot and link has a diagram that realizes the sequence. We sometimes restrict to sequences that are universal just for knots but not for multi-component links.
\end{defn}

\begin{figure}[h]
\begin{center}
\includegraphics[scale=0.75]{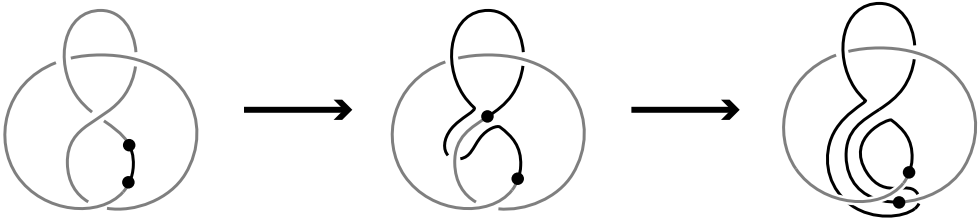}
\caption{A move similar to \cite[Figure 3]{Oz}}\label{ozawa-new}
\end{center}
\end{figure}

This paper is an investigation into which sequences are universal for knots and/or links. This research  was motivated by a fact from \cite[Remark 2.2 and Figure 3]{Oz}. (See \figurename~\ref{ozawa-new}, which is a move similar to \cite[Figure 3]{Oz}.)

\begin{thm}\label{thm:two-coloring} Every knot has a projection that can be decomposed into two sub-arcs such that each sub-arc never crosses itself.
\end{thm}

To see why the theorem holds, start with a short black arc that avoids the crossings and make the rest of the projection grey.  Then  extend the  black arc by moving one of its vertices along the knot projection, as in \figurename~\ref{ozawa-new}. If we come to a crossing and it is not yet colored, extend the black through the crossing. If we come to a crossing and the other strand is already black, push this strand along ahead of the black coloring as we extend the black coloring. Continue in this manner until we have a valid 2-coloring. 

 This theorem easily implies that every knot has a projection with at most four odd-sided regions. Take the projection provided by the theorem. Let the two points where the knot changes color be called {\it transition vertices}. Any face that does not contain a transition vertex on one of its edges must have an even number of edges that alternate between the two colors as we travel around the boundary of the face.  A face that contains one transition vertex on an edge must have an odd number of edges. Since there are two transition vertices and each can appear on an edge incident to two faces, there are at most four odd-sided faces. In Theorem \ref{thm:twoodd}, we improve this to show that in fact every knot has a projection with exactly two odd-sided faces.

 This result  causes one to consider what the restrictions are on the values of $n$ that occur for the faces in  projections of knots and links.

Let us begin with some elementary observations.  Given a knot diagram $P$, let $p_i$ be the number of faces of $i$ sides. Note that $p_1 = 0$. Considering the previously mentioned 4-valent graph on the 2-sphere, we have: 

\begin{equation} e = \frac{2p_2+3p_3+ 4p_4+ \dots}{2}
\end{equation}
 
\begin{equation} f= p_2+p_3+p_4 +\dots
\end{equation}

\begin{equation} e=4v/2
\end{equation}

By Euler Characteristic, we can plug this into $v-e+f=2$ to obtain:

\begin{equation}{\label{eqn:Euler}}
2p_2 + p_3= 8 + p_5+ 2p_6 + 3p_7 + \dots
\end{equation}

In graph theory and geometry, there is a long history, dating back to \cite{Eberhard}, of investigations into those sequences of values that satisfy Equation \ref{eqn:Euler} and that represent actual $4$-valent graphs on the sphere. However, graph theorists and geometers, who are interested in these graphs as the 1-skeletons of convex polytopes,  typically assume that $p_2= 0$, which we will not usually want to do. 

For instance, in \cite{Grunbaum}, Gr\"unbaum proved that for any collection of integers with $p_2$ = 0 that satisfies Equation \ref{eqn:Euler}, there exists a choice of $p_4$ such that there is a planar 3-connected $4$-valent graph realizing these values. This theorem is known as Eberhard's Theorem. In \cite{Jeong}, this result was extended to show that there is a value for $p_4$ so that the resulting graph is the projection of a knot, rather than a link.

For our purposes, Equation \ref{eqn:Euler} provides several useful pieces of information. First, for any diagram of a knot or link, either $p_2$ or $p_3$ is nonzero. In other words, any sequence that is realized by a nontrivial reduced projection must begin with  a 2 or a 3. Second, the number of $4$-gons is not constrained by Euler characteristic. Third, by considering this equation mod 2, the number of odd-sided regions must always be even.

In Section 2, we begin by proving the result mentioned above, that every knot has a diagram with exactly two odd-sided faces. Moreover, we show that the odd-sided regions can be chosen to be triangles.

Putting this result in terms of what it implies in the language introduced above, the sequence $(2, 3, 4, 6, 8, \dots)$ is universal for knots. We further obtain upper bounds on the least number of odd regions that can occur in a projection of an $n$-component link.

In the opposite direction, we prove that the integers in a  sequence realized for any knot cannot have a nontrivial common divisor. For example, it can never be the case that all the faces are even-sided. This is not true for links. In fact, a 2-component link is realized in a projection with all even-sided faces if and only if both components are individually trivial.

Note that because every universal sequence must begin with a 2 or a 3, this theorem only excludes the two possibilities for universal sequences of knots that all the integers are even and that all the integers are divisible by 3.

In Section 3, we show that certain sequences are universal. For instance, we prove that $(3,5,7,\dots)$ is universal and that both  $(2, n, n+1, n+2\dots)$ for any integer $n \ge 3$ and $(3,n,n+1,n+2, \dots)$ for any integer $n \ge 4$ are universal.

In fact, universal sequences need not be infinite. We prove that the sequence $(3,4,n)$ is universal for all $n \ge 5$.  Note that such diagrams have no bigons and as such, are called {\it lune-free diagrams} in \cite{EHK}. In that paper, the authors produce infinitely many lune-free diagrams, and show that such diagrams are realized for any number of vertices greater than 7. It is somewhat surprising that lune-free diagrams, particularly ones as simple as for instance $(3,4,5$)-diagrams, yield projections for all knots.

Finally, we prove that the sequence $(2,4,5)$ is universal. We end with some open questions.

Note that all of the proofs in this paper that sequences are universal do not alter the crossings in the original projection. Because of this, the results appearing here can be stated in exactly the same form for any of the categories of virtual, flat, welded, and singular knots. 

Throughout this paper, we use the fact that any projection can be made a projection of the trivial knot by changing crossings, a proof of which appears in several elementary knot theory texts.

The result from \cite{Oz} mentioned above together with results from \cite{LJ} which we utilize in the next section are both generalized to spatial embeddings of graphs in \cite{S}.

\section{Minimizing Odd Regions}

\begin{thm}{\label{thm:twoodd}}Every knot has a diagram with exactly 
two odd-sided faces, which can be made to be triangles.
\end{thm}

\begin{proof} 
We take a $2$-colored diagram $D$ of a given knot $K$ such that neither colored arc crosses itself
as we mentioned in Section~\ref{S:intro}. (See \figurename~\ref{ozawa-new}.) 
All faces must have an even number of edges except 
for those that have transition vertices on their edges. 
Then the diagram $D$ has two or four odd-sided regions 
as we discussed in Section~\ref{S:intro}. 
If $D$ has four odd-sided regions, we apply the following 
procedure.

Take a parallel copy $D'$ of $D$ such that 
newly created crossings have the following crossing information. 
\begin{itemize} 
\item[(1)] 
For a crossing between $D$ and $D'$, 
the edge of $D'$ is lower than that of $D$, and 
\item[(2)]
For a crossing both of whose edges belong to $D'$, 
the gray edge is lower than the black edge. 
\end{itemize}
The diagram $D \cup D'$ represents a split two component 
link with $K$ and a trivial knot, and
it is also $2$-colored in a natural way. 
(See \figurename~\ref{two-odd}.) 
Taking a band $b$ near one of transition vertices
and performing a band surgery along $b$ as in \figurename~\ref{two-odd}, 
we obtain a $2$-colored diagram of $K$ such that there exists a region 
whose boundary edges contain both of two vertices. 
This is because this band surgery induces
the connected sum of $K$ and the trivial knot.
Then we have a diagram of $K$ with exactly two odd-sided regions.
By applying moves as in \figurename~ \ref{triangle},
we make these two odd-sided regions into two $3$-sided regions and even-sided regions.

\begin{figure}[h]
\begin{center}
\includegraphics[scale=0.7]{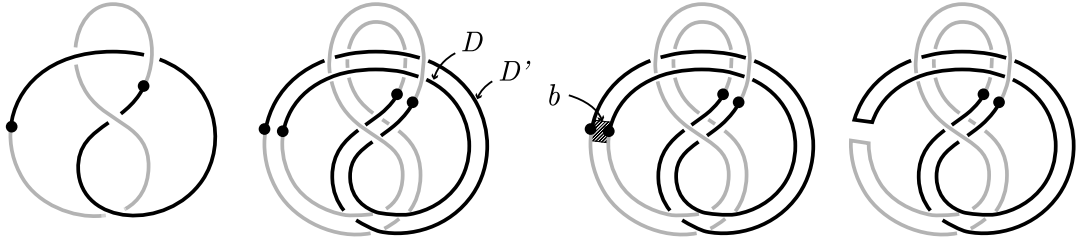}
\caption{Example for $4_1$}
\label{two-odd}
\end{center}
\end{figure}

\begin{figure}[h]
\begin{center}
\includegraphics[scale=0.7]{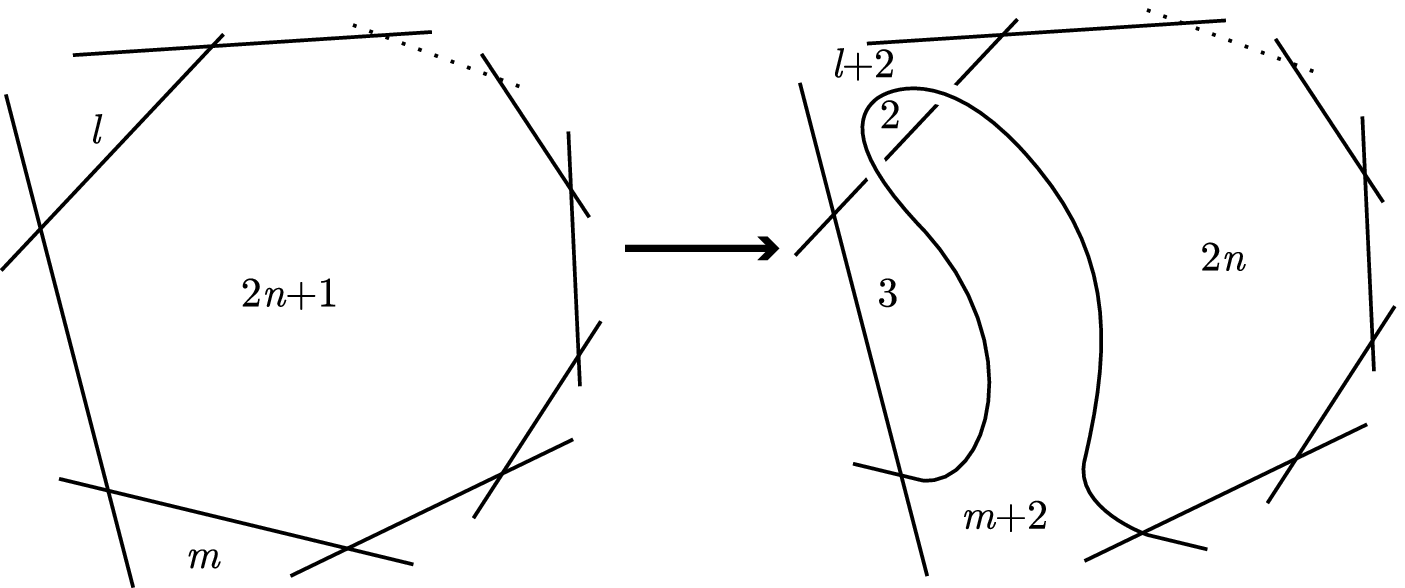}
\caption{}
\label{triangle}
\end{center}
\end{figure}
\end{proof}

\begin{thm}{\label{thm:nodd}}Every link $L$ of $n$ components has a diagram with at most $n$ odd-sided faces if $n$ is even and $n+1$ odd-sided faces if $n$ is odd. The odd-sided faces can be made to be triangles.
\end{thm}

\begin{proof} Start with any projection of the $n$-component link $L$. Then we can rearrange it to construct a region $R$ that has $n$ edges, each edge coming from a different component. Put a transition vertex at the center of each edge, so that as we move clockwise around the boundary of the region, each edge is first black and then gray. So at least on the boundary of this region, there are no crossings of a single color. Now we will rearrange the projection so that all crossings have only one color. As in the proof of Theorem \ref{thm:two-coloring}, start with one of the components and extend its black section  until we have a valid 2-coloring of this component, pushing strands ahead as necessary.  Repeat this process with each link component in turn. We obtain a 2-coloration  of the link such that  no colored edge crosses itself. The $n$-gon $R$ remains and it contains a transition vertex between colors on each of its edges.

As in the proof of Theorem \ref{thm:twoodd}, all faces must have an even number of edges except for those that have transition vertices on their edges. If there are an odd number of edges with transition vertices on them, the face must have an odd number of edges. If there are an even number of edges with transition vertices, the face must have an even number of edges.

Applying the doubling procedure as in the proof of Theorem \ref{thm:twoodd}, the second transitional vertex for each of the $n$ components can be made to share a quadrilateral with the first transitional vertex for that component. Note that $R$ also survives this procedure. Then each of these second transitional vertices can generate at most one odd-sided face, and $R$ itself is odd-sided if and only if $n$ is odd. Hence the theorem results. As in the previous proof, we can convert the odd-sided regions into triangles.
 \end{proof}

\begin{thm} A sequence with a nontrivial common divisor $k$ can only realize links with at least $k$ components.
\end{thm}

\begin{proof} Suppose that $k > 1$ divides all of the integers in the sequence. Choosing any region $R$ to begin, its number of edges $n$ must satisfy $n=kr$ for some integer $r$. Starting on one edge and traveling around the boundary of $R$, we  label each edge with consecutive integers $1,2,\dots, k$, and then repeat $r$ times as we cycle around the boundary of $R$. For each edge that has been labelled, we can extend its label to the edge that is opposite it at each of its endpoints. In this way, we obtain labelings on two of the edges of each regions that shares just a vertex with $R$ and three of the edges of a region that shares an edge with $R$. Extend these labels to an adjacent region so that the union of labeled regions remains topologically a disk. Repeat this process until the entire diagram is labeled. 

To see that this method of labeling the edges in the diagram is well-defined, suppose that we have labeled it appropriately for $n$ faces, the union of which is the disk $D$, and we are adding in an additional face $F$. The face $F$ intersects the boundary of $D$ in an arc $\alpha$ made up of the union of a collection of edges. At a vertex $v$ internal to the arc $\alpha$, the labeling of the face in $D$ opposite to $F$ at $v$ forces the two edges on the boundary of $D$ at $v$ to have consecutive labels. The labelings of the faces in $D$ containing those two edges force the labelings on the adjacent edges in $\alpha$ to be consistent. At a vertex $v$ that is an endpoint of $\alpha$, the labeling on the face that contains the edge in $\alpha$ that contains $v$ forces a consistent labeling of the subsequent edge in $F$ that can be extended around $F$ to obtain a consistent labeling of $D \cup F$, as desired.

Hence, all the edges that make up a given component receive the same label. Since $k$ labels appear, the diagram must have at least $k$ components. 
\end{proof} 

Thus, in particular, knots cannot have a diagram realized by an even sequence. The next lemma determines which  links can be so realized.

\begin{lemma} An $n$-component link $L$ is realized by an even sequence if and only if its components are all trivial and the components of $L$ can be subdivided into two nontrivial sets such that the complement of either set is a trivial link.
\end{lemma}

\begin{proof} If an $n$-component link is realized by an even sequence, then as in the proof of the preceding theorem, we can label the diagram with two labels such that each component receives a consistent label. Since crossings only occur between distinctly labeled components, no component crosses itself and hence, all components are trivial. The two labels divide the components into two sets, such that the complement of either set has no crossings and therefore is the trivial link.

 In the other direction, Theorem 2 of \cite{LJ} implies that  there exists a single projection of $L$ such that each of the two subsets of components appears with no self-crossings. This implies that the two sets must alternate on the edges  as we travel around the boundary of any region and therefore all of the regions in that projection must be even-sided.
\end{proof}

\section{Universal Sequences}
 
\begin{thm}The following sequences are universal for knots and links.
\begin{enumerate}
\item[i)] 
$(3,5,7,9, \ldots)$
\item[ii)] 
$(3,n,n+1,n+2, \dots)$ for any integer $n \ge 4$.
\item[iii)] 
$(2, n, n+1, n+2\dots)$ for any integer $n \ge 3$.
\end{enumerate}
\end{thm}

\begin{proof}
i)
Let $T$ be the diagram of the trivial knot as in \figurename
 \ \ref{triv} with $k=1$.
Any even sided region of any diagram $D$ can be changed 
into ten odd sided regions by attaching $T$ to a region $R$ 
as in \figurename~\ref{(3,n,n+1)} without increasing the number of 
even sided regions. 
We apply this deformation to each even sided regions of $D$.
It is easy to check the resulting diagram is also a diagram of $L$.

\begin{figure}[htbp]
\begin{center}
\includegraphics[scale=0.8]{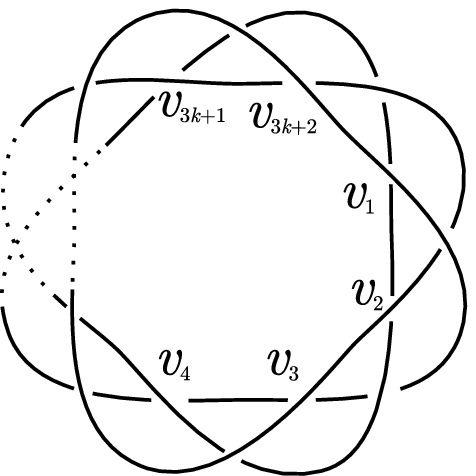}
\caption{}
\label{triv}
\end{center}
\end{figure}

\begin{figure}[htbp]
\begin{center}
\includegraphics[scale=0.6]{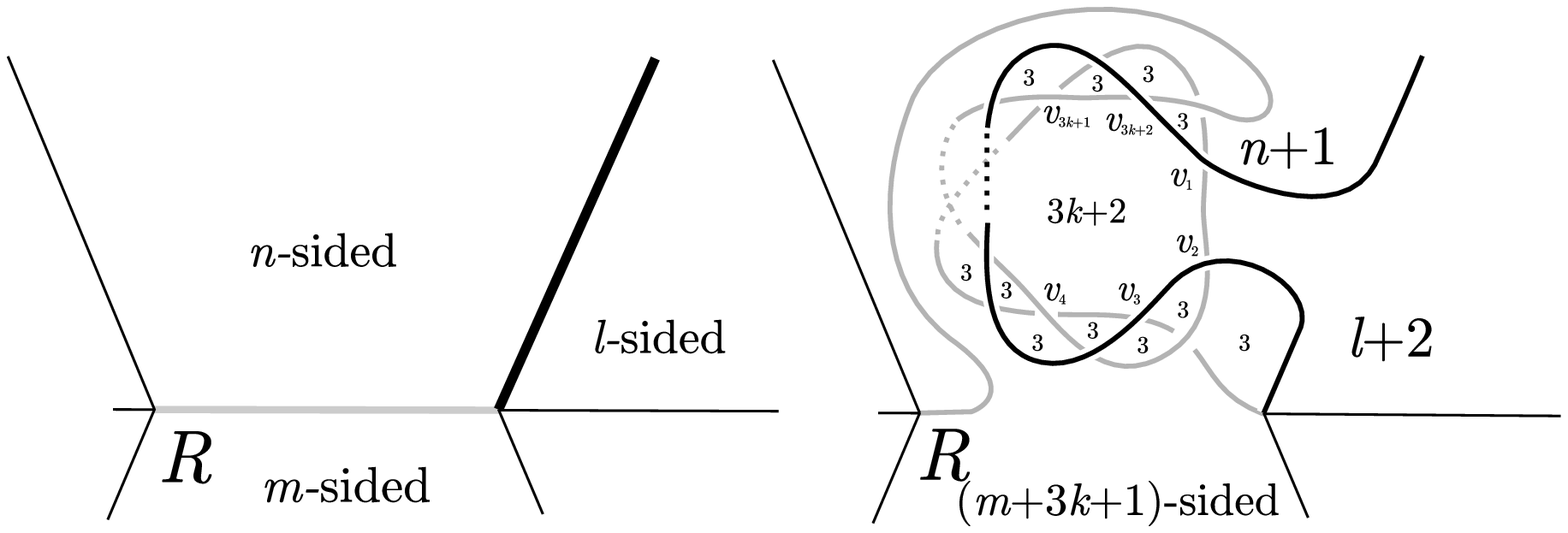}
\caption{}
\label{(3,n,n+1)}
\end{center}
\end{figure}
 
ii) 
Let $k$ be a natural number such that $3k+2 \ge n$.
Let $D$ be a diagram of a knot and $T$ the diagram of the trivial knot   
with $(6k+4)$ crossings as in \figurename~\ref{triv}.
By attaching $T$ to a region $R$ of $D$ as in \figurename~\ref{(3,n,n+1)}, 
$(6k+2)$ triangles are newly created 
and the number of edges of $R$ increases by $3k+1$ 
without decreasing that of each region of $D$. 
The desired diagram can be obtained by applying this operation 
to each region the number of whose edge is less than $3$.

iii) 
If there exist a region $R$ with fewer than $n$ edges, we apply the move as in \figurename~\ref{(2,n,n+1)}.
This move creates bigons and increases the number of edges of $R$ 
without decreasing that of any other region of $D$.

\begin{figure}[htbp]
\begin{center}
\includegraphics[scale=0.8]{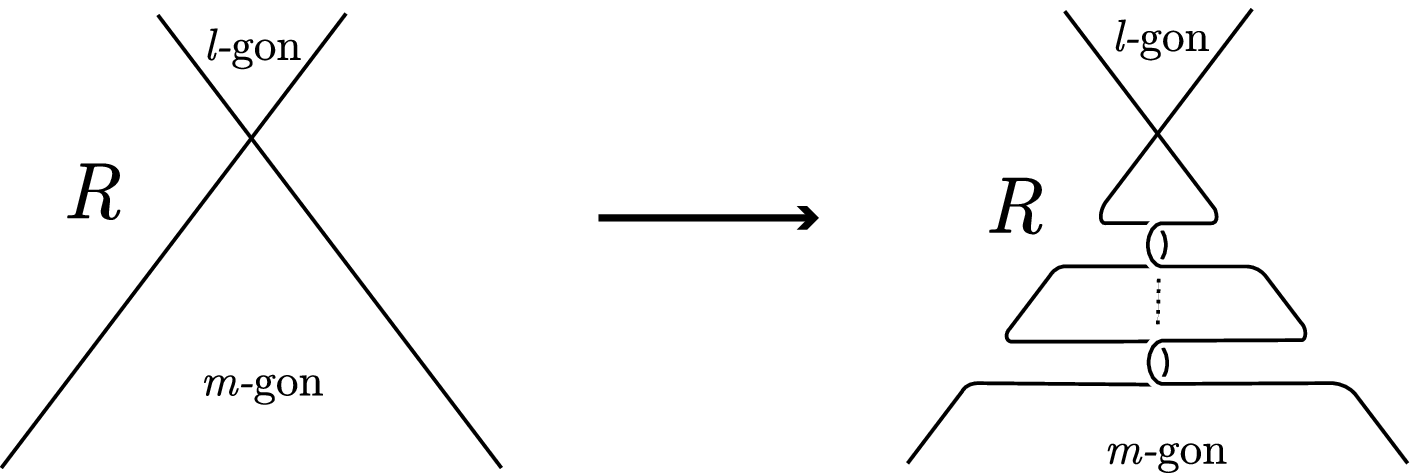}
\caption{}
\label{(2,n,n+1)}
\end{center}
\end{figure}

\end{proof}

\begin{thm} The sequence $(3,4,5)$ is universal for knots and links.
\end{thm}

\begin{proof} Consider the projection $Q$ in Figure \ref{fig:(3,4)-proj}. It is made up of eight triangles and a collection of quadrilaterals. Note that by continuing to spiral in the two possible directions, such a projection could be created with eight triangles and a grid of quadrilaterals at the center that is arbitrarily large in both the horizontal and vertical directions. For convenience, we continue to call such a projection $Q$, regardless of the number of quadrilaterals. (Note: In \cite{Enns}, it was proved that for any $m >1$, there exists a $4$-valent graph with eight triangles and $m$ quadrilaterals. However, it is not necessarily a knot projection.)
 
\begin{figure}[htbp]
\begin{center}
\includegraphics*[height=4.3in]{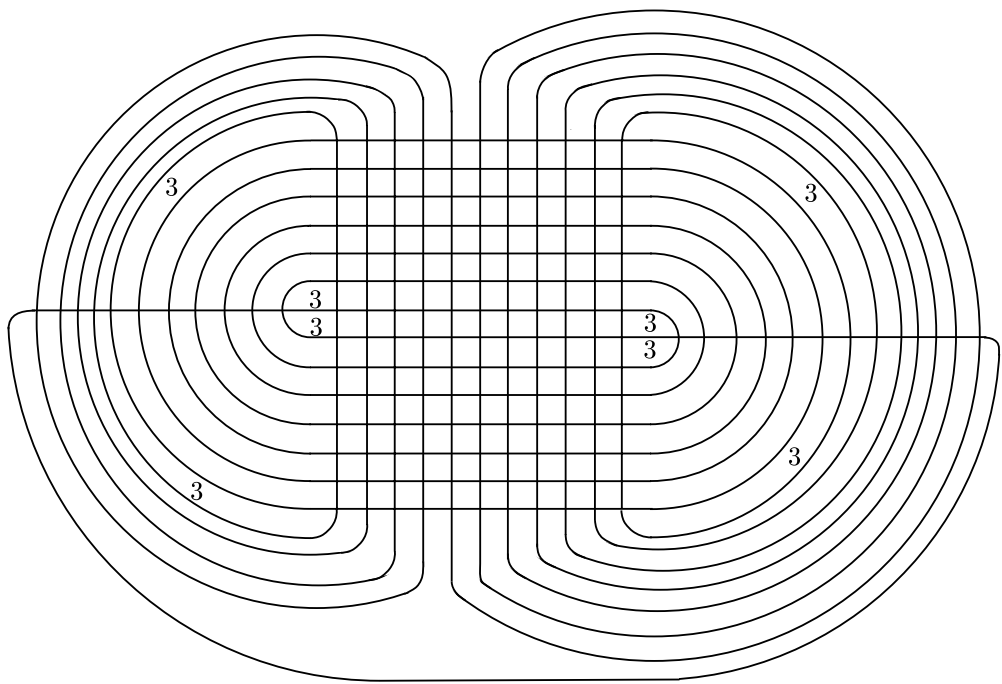}
\end{center}
\caption{\label{fig:(3,4)-proj}A $(3,4)$-projection on the sphere.}
\end{figure} 

Let  $P$ be a reduced projection of any given knot or link $K$. Isotope the projection $P$ to a projection $P'$  that  lies in the grid of horizontal and vertical lines passing through the integers on the x and y axes in the plane. All crossings will occur at integer lattice points.  The transition point in $P'$ between a vertical line segment in the grid and a horizontal line segment in the grid is called a {\it corner}. Place over $P'$ a projection of the trivial knot in the shape of $Q$ so that the vertices of the quadrilateral grid at the center of $Q$ have half-integer coordinates and so that $Q$ has  enough squares in both directions of its central quadrilateral grid to cover the entire projection $P$. Assume that the central grid of $Q$ is fine enough so that if $P'$ intersects a given quadrilateral $S$, it does so either as a single arc intersecting opposite edges of $S$, or in a single arc containing one corner that intersects two adjacent edges of $S$, or it intersects $S$ in a crossing, with the four arcs coming out of the crossing intersecting the four edges of $S$. (See 
%the upper half of 
\figurename~\ref{fig:(3,4,5)-1(1)}.) We also assume that the grid is fine enough so that there is at least one quadrilateral that does intersect $P'$ in a single arc intersecting opposite edges and that one of quadrilaterals adjacent to it does not intersect $P'$ as shown in the 
%bottom 
left of \figurename~\ref{fig:(3,4,5)-1(2)}. 

%\begin{figure}[htbp]
%\begin{center}
%\includegraphics*[height=3in]{figures/8.eps}
%\end{center}
%\caption{\label{fig:(3,4,5)-1}The part where $P'$ and $Q$ are composed.}
%\end{figure} 

\begin{figure}[htbp]
\begin{center}
\includegraphics*[height=1.2in]{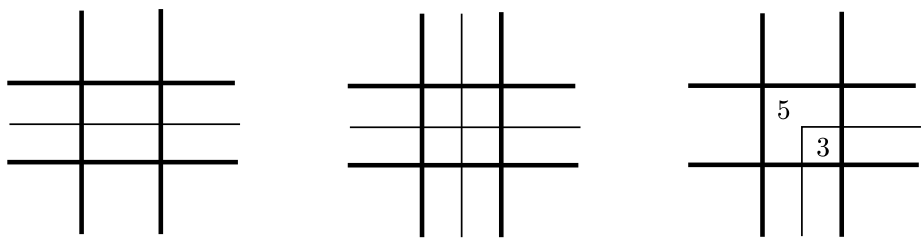}
\end{center}
\caption{\label{fig:(3,4,5)-1(1)}The three situations where 
$P'$ intersects a quadrilateral of $Q$.}
\end{figure} 

\begin{figure}[htbp]
\begin{center}
\includegraphics*[height=1.5in]{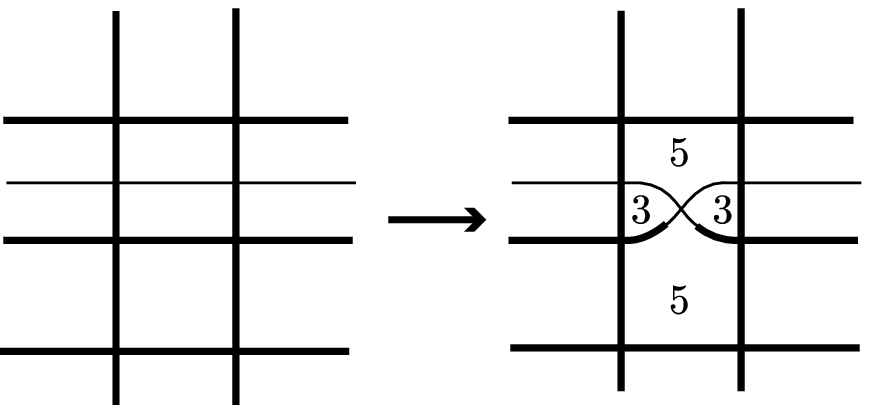}
\end{center}
\caption{\label{fig:(3,4,5)-1(2)}The part where 
$P'$ and $Q$ are composed.}
\end{figure} 

The link that results by including $Q$ with $K$ is thus in a $(3,4,5)$-diagram. Now, to compose the grid with one of the components of the link, we take an arc of $P'$  that is passing through a quadrilateral, intersecting opposite edges, and we add a crossing between the arc from $P'$  and one of the edges from $Q$  that it does not intersect. This replaces three quadrilaterals with  two triangles and two pentagons. (See the 
%bottom 
right of \figurename~\ref{fig:(3,4,5)-1(2)}.) The resulting link is  the composition of the link $K$  and the trivial knot and is therefore  the link $K$ in a (3,4,5)-projection.
\end{proof}

\begin{thm} The sequence $(3,4,n)$ is universal for knots  for all $n \ge 5$. 
\end{thm}

\begin{proof} To prove that there exists a $(3,4,n)$-diagram for $K$ in the case $K$ is a knot, begin with the same projection $P'$ as appeared in the previous proof overlayed by the same $Q$. Note that $P'$ has an even number of corners. Take an additional $n-5$ copies of $P'$, each with crossings chosen to yield a trivial knot and each laid one on top the other such that at the corners of $P'$, they cross one another as in Figure \ref{fig:corner}, which depicts the case $n=8$. 

\begin{figure}[htbp]
\begin{center}
\includegraphics*[height=2in]{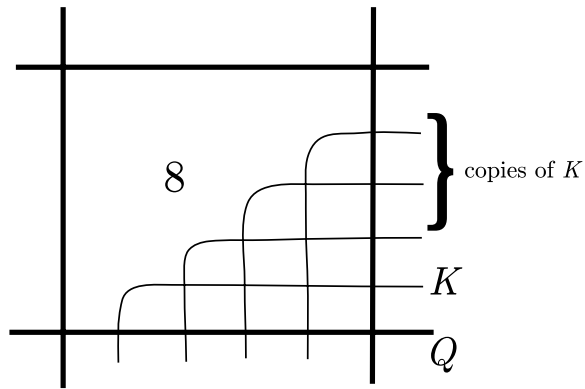}
\end{center}
\caption{\label{fig:corner}The new projection at a corner.}
\end{figure} 

  The resulting link projection is a $(3,4,n)$-link projection. Finally, to connect the projection of $K$ with the grid and with the $n-5$ new trivial components, we choose a quadrilateral of $Q$ through which the $n-4$ components pass straight through, and we put in a twist as in Figure \ref{fig:twist}. 

\begin{figure}[htbp]
\begin{center}
\includegraphics*[height=2.5in]{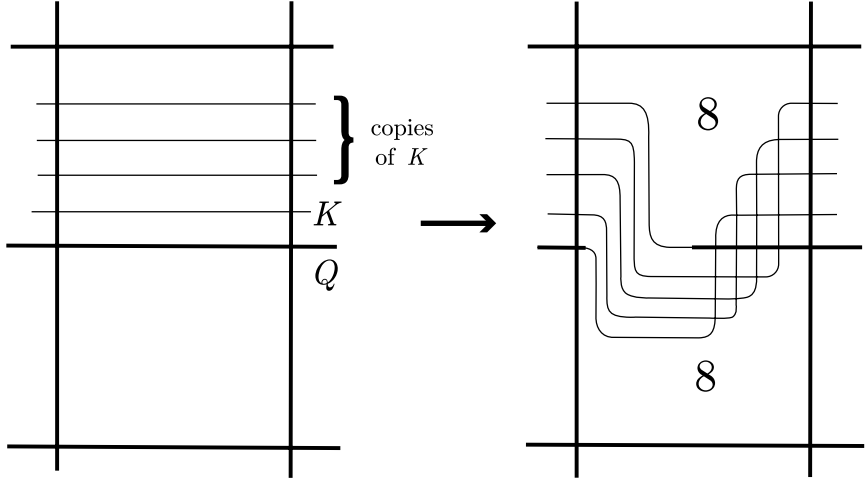}
\end{center}
\caption{\label{fig:twist}Twisting to compose.}
\end{figure} 

The resulting knot is in fact $K$ in a $(3,4,n)$-projection.
\end{proof}

Finally, we prove:

\begin{thm} The sequence $(2,4,5)$ is universal for knots and links.
\end{thm} 

\begin{proof} 
Start with the trivial knot projection $Q$ 
in \figurename~\ref{fig:(3,4)-proj},
that is entirely made up of eight $3$-gons and a collection of $4$-gons.
By deforming $Q$ near $3$-gons as 
in the upper and bottom of \figurename~\ref{fig:(2,4,5)proj}, 
we turn $Q$ into a projection 
that is entirely made up of $2$, $4$, and $5$-gons. 
(We also call this projection $Q$.)

For a given projection of a knot or link  $K$, we deform it into 
a lattice projection $P$ and take a second copy $P'$ 
of the projection $P$, but chose the crossings for P' to result in a trivial link sitting on top of the projection $P$.
By connecting $P$ and $P'$ as in \figurename~\ref{fig:(2,4,5)parallel}, 
we obtain a projection $R$ of $K$. Note that if $K$ has $n$ components, we connect $P'$ to $P$ at $n$ places, one for each component.
We then set $Q$ on the top of $R$ 
so that the following conditions are satisfied.

\begin{figure}[htbp]
\begin{center}
\includegraphics*[height=3.5in]{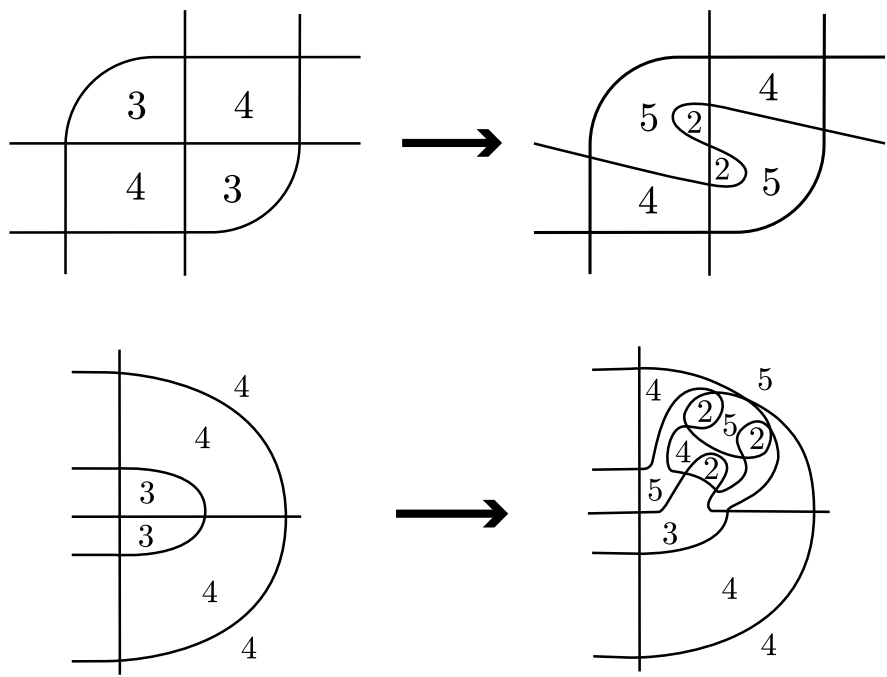}
\end{center}
\caption{\label{fig:(2,4,5)proj}From a $(3,4)$-projection to a $(2,4,5)$-projection.}
\end{figure}

\begin{figure}
\begin{center}
\includegraphics*[scale=0.7]{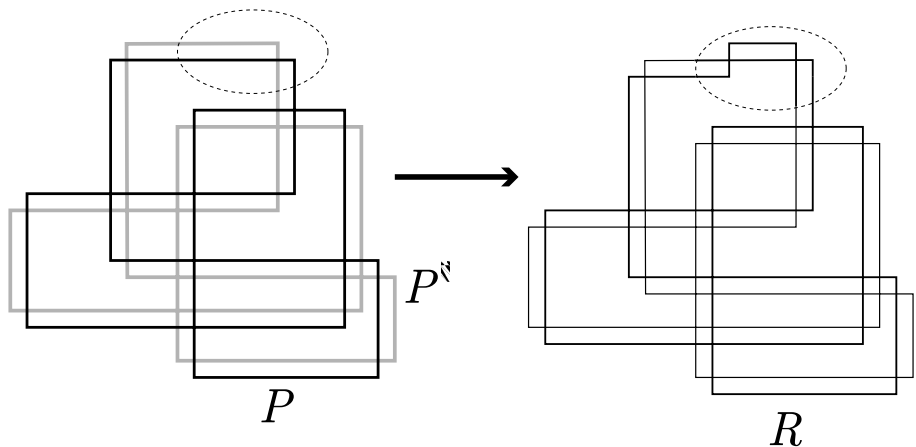}
\end{center}
\caption{\label{fig:(2,4,5)parallel}The diagram $R$ obtained from $P$ and $P'$.}\end{figure}

\begin{itemize}
\item[(1)] 
$R$ is entirely covered by a portion of $Q$ tiled by quadrilaterals 
as in \figurename~\ref{fig:(2,4,5)onR}. 
(At a  neigbourhood of each point of a vertical line of $P$, a part of $R$ 
is contained in one quadrilateral, and at a  neigbourhood
of each point of a horizontal line of $P$, 
a part of $R$ is contained in two quadrilaterals.)

\item[(2)]
Each corner of $R$ is contained in  two quadrilaterals of $Q$.
(See the left of \figurename~\ref{fig:(2,4,5)corner}.)
The lower quadrilateral is cut into two triangles, a quadrilateral and
a pentagon, and the upper one into a triangle and a pentagon.

\item[(3)] 
The parts where $P$ and $P'$ are connected are each contained in three
quadrilaterals of $Q$. 
(See the left of \figurename~\ref{fig:(2,4,5)compose}.) 
Each of upper and lower quadrilaterals is cut into 
a triangle and a pentagon, and the middle one is cut into four quadrilaterals.
\end{itemize}

\begin{figure}[htbp]
\begin{center}
\includegraphics*[height=3.3in]{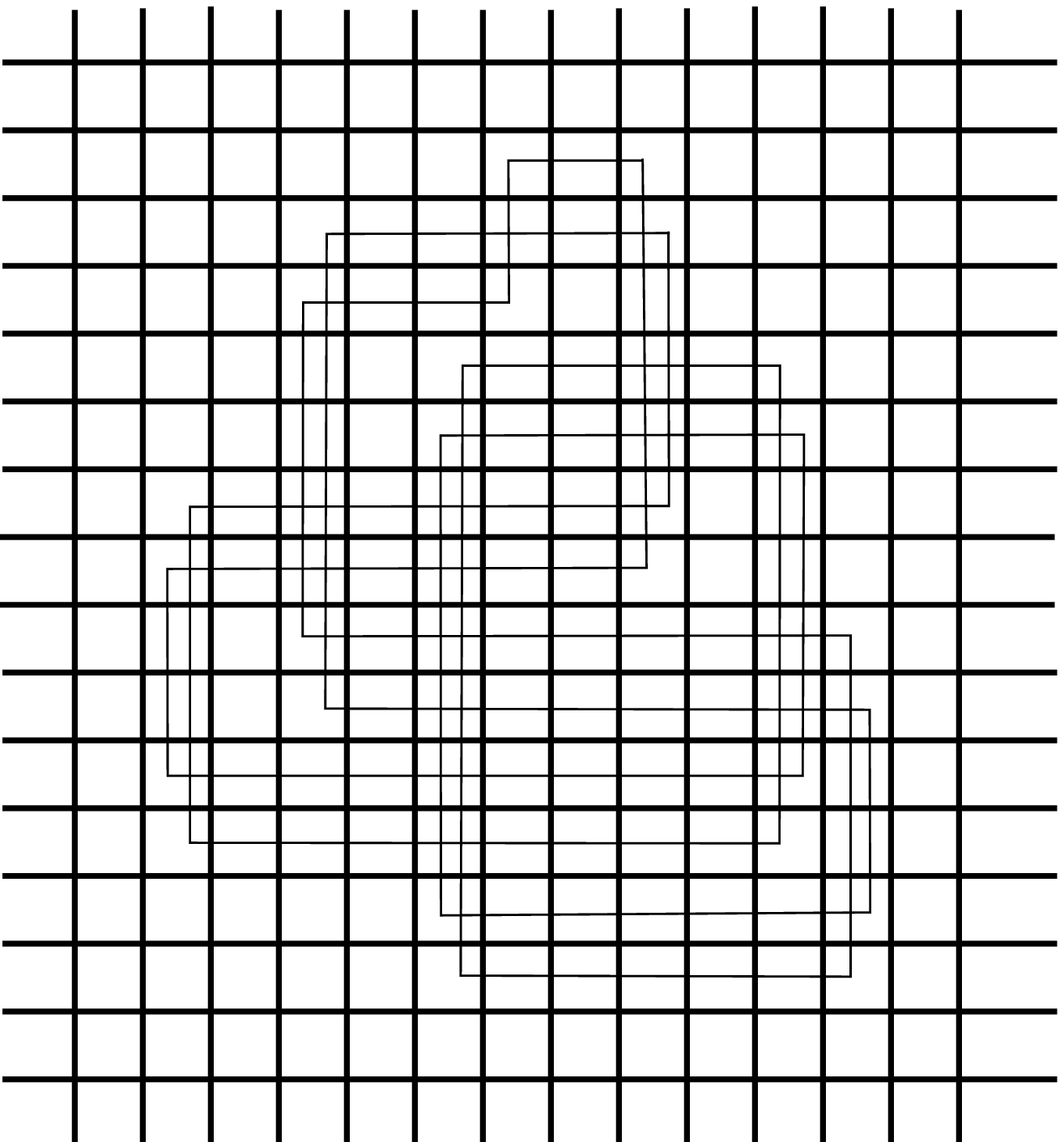}
\end{center}
\caption{\label{fig:(2,4,5)onR}Setting $Q$ on the top of $R$.}
\end{figure}

Since the projection $R$ and $Q$ is entirely made up of 
$2$, $3$, $4$ and $5$-gons, we need to erase all $3$-gons. 
At each corner we apply the move as in \figurename~\ref{fig:(2,4,5)corner},
and at the parts where $P$ and $P'$ are connected 
we apply the move as in \figurename~\ref{fig:(2,4,5)compose}.
Then we obtain a projection $R$ and $Q$ with only $2$, $4$ and $5$-gons.
The final step is to connect $R$ to $Q$ without increasing $3$-gons.
By altering the projection as in \figurename~\ref{fig:(2,4,5)RandQ}, 
we obtain a $(2,4,5)$-projection of K.

\begin{figure}[htbp]
\begin{center}
\includegraphics*[height=1.8in]{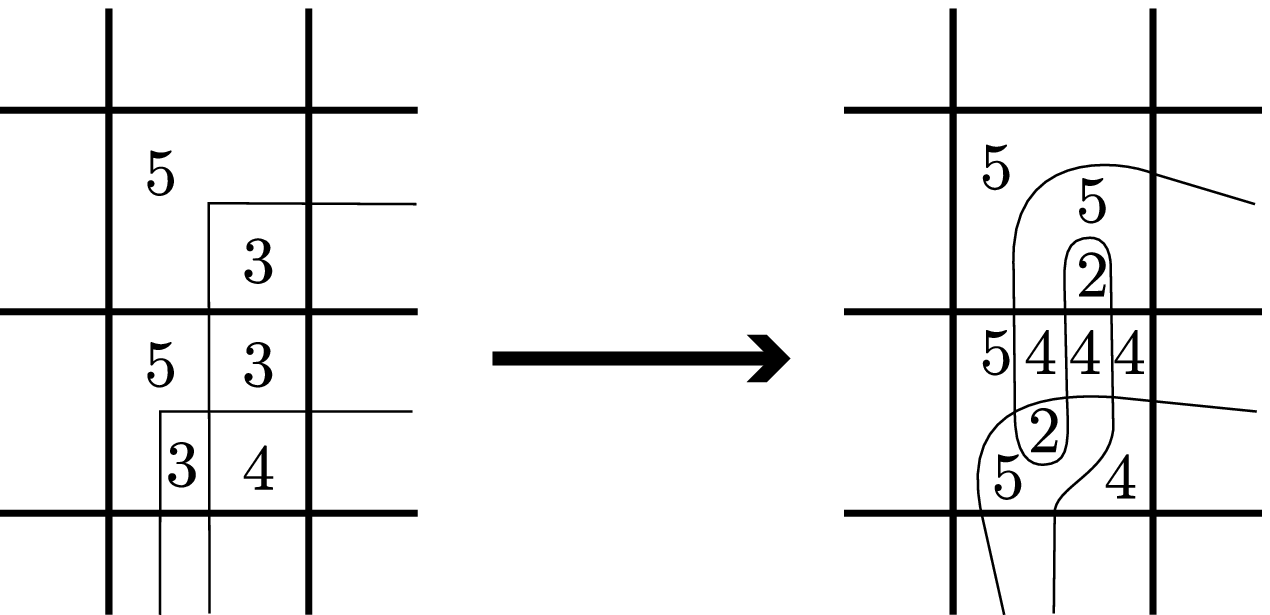}
\end{center}
\caption{\label{fig:(2,4,5)corner}The new projection at a corner.}
\end{figure} 

\begin{figure}[htbp]
\begin{center}
\includegraphics*[height=2.4in]{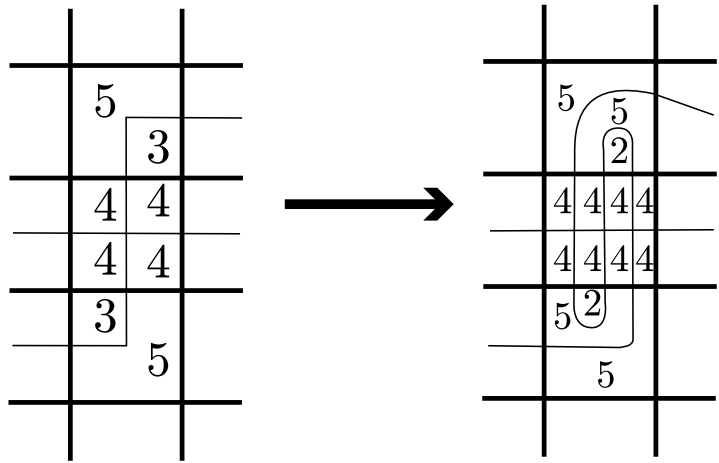}
\end{center}
\caption{\label{fig:(2,4,5)compose}The part where the copies of $K$ are composed.}
\end{figure} 

\begin{figure}[htbp]
\begin{center}
\includegraphics*[height=1.5in]{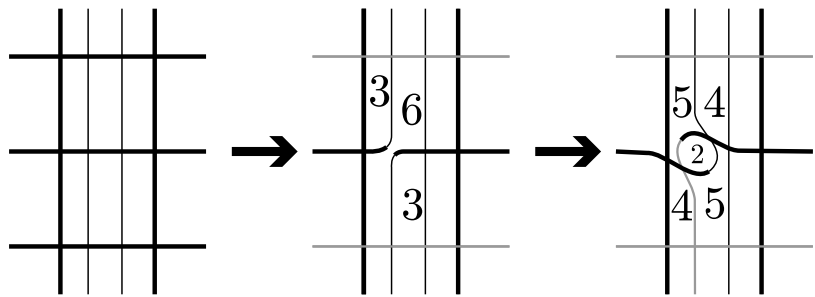}
\caption{\label{fig:(2,4,5)RandQ}
The part where $R$ and $Q$ are composed.}
\end{center}
\end{figure}

\end{proof}

\medskip

{\bf Questions:}

\begin{enumerate}

\item Does there exist a universal sequence for knots or for links of just two integers? Note that it would have to be of the form $(2,n)$ where $n$ is odd or $(3,m)$ where $m$ is not divisible by 3. 

\item In particular, one would like to know whether or not (3,4) is a universal sequence. Note that for a (3,4)-projection, there must be exactly 8 triangles, however as in Figure \ref{fig:(3,4)-proj}, the number of quadrilaterals in such projections can be arbitrarily large. Since the crossings can be chosen to be alternating, and since \cite{K}, \cite{M},\cite{T} proves that as the crossing number goes up the resulting knots will all be distinct, we know that the resulting knots form an infinite collection. So there are infinitely many knots realizing (3,4)-projections. But does this infinite collection contain all knots?

\item Can one prove that there exists links for which Theorem \ref{thm:nodd} is best possible?

\end{enumerate}

\end{document}